\documentclass[12pt,letterpaper, reqno]{amsart}
\usepackage{amssymb,latexsym, amsmath, amsxtra}
\usepackage{graphicx}
\usepackage{array}

\newtheorem{theorem}{Theorem}[section]
\theoremstyle{plain}

\theoremstyle{definition}

\theoremstyle{remark}

\numberwithin{equation}{section}
\numberwithin{theorem}{section}
\numberwithin{table}{section}
\numberwithin{figure}{section}

\newcommand{\C}{\mathbb C}

\newcommand{\R}{\mathbb R}
\newcommand{\Z}{\mathbb Z}

\newcommand{\cj}[1]{\overline{#1}}
\newcommand{\abs}[1]{\left| #1 \right|}
\newcommand{\cbr}[1]{\left\{ #1 \right\}}
\newcommand{\step}[2]{\medskip\noindent{\bf Step #1. \enspace}{\emph{#2}}}

\newcommand{\RHS}{\mbox{RHS}}
\newcommand{\AppFE}{\mbox{AppFE}}

\def\({\left(}
\def\){\right)}

\begin{document}
\title{Maass forms on GL(3) and GL(4)}
\author{David W. Farmer, Sally Koutsoliotas, and Stefan Lemurell}



\begin{abstract}
We describe a practical method for finding an $L-$function
without first finding the associated underlying object.
The procedure
involves using the Euler product and the approximate
functional equation in a new way.  No use is made of the functional
equation of twists of the $L-$function.
The method is used to find a large number of
Maass forms on $SL(3,\mathbb Z)$ and to give the first
examples of Maass forms of higher level on $GL(3)$, and on
$GL(4)$ and $Sp(4)$.

\end{abstract}

\address{
{\parskip 0pt
American Institute of Mathematics\endgraf
farmer@aimath.org\endgraf
\null
Bucknell University\endgraf
koutslts@bucknell.edu\endgraf
\null
Chalmers University of Technology and\endgraf
University of Gothenburg\endgraf
sj@chalmers.se\endgraf
}
  }

\thanks{The authors thank Ce Bian, Andrew Booker, Brian Conrey,
Stephen D.~Miller,  Michael Rubinstein, 
Ralf Schmidt, and Nicolas Templier
for helpful discussions.}

\maketitle

\section{Introduction}

We describe a general approach to finding an $L-$function given only a limited
amount of information about its functional equation
and Euler product.  We illustrate the method by
locating $L-$functions associated to Maass forms on $SL(3,\Z)$ and some of its
subgroups,  and on $Sp(4,\Z)$ and $SL(4,\Z)$.

\subsection{Summary of the method and results}

We describe an approach to solving the following problem:
find all $L-$functions satisfying a given functional equation.
The precise assumptions, which also involve an Euler product and a
Ramanujan bound on the coefficients, are given  in
Section~\ref{sec:axioms}.  Our results show that, at least for
$L-$functions of 
degree $d\le 4$, this problem can be solved in a practical way,
without an assumption of a functional equation for twists of the
$L-$function and without first finding an arithmetic object that 
gives rise to the $L$-function.

Our approach involves a new way 
of extracting information from
the smoothed approximate
functional equation for an $L-$function.  This formula, given in
Theorem~\ref{thm:formula}, is a standard tool in analytic number
theory.  The new ingredient makes use of the fact that a variety of test
functions can appear in the formula.  
This allows us to create a relatively small system of equations
which the Dirichlet coefficients of the
unknown
$L-$function must satisfy.  Solving that system provides the
missing information about the $L-$functions with the given
functional equation.

Using this method, we have found more than 2000 $L-$functions of Maass
forms on $SL(3,\Z)$ and several dozen $L-$functions of Maass forms
on congruence subgroups of $SL(3,\Z)$.  The method also shows
the non-existence of small eigenvalues for Maass forms on
$SL(3,\Z)$ in a manner that is more effective than previous
results using the explicit formula~\cite{Mil1}.

We have also found more than 200 $L-$functions of Maass forms
on $SL(4,\Z)$ and several dozen for $Sp(4,\Z)$.  The smallest
$Sp(4,\Z)$ example has the surprising property that its
first zero on the critical line has a larger imaginary
part than the first zero of the Riemann zeta-function~\cite{Bo}.
We have also used this method to find, with proof, the smallest
conductor of a hyperelliptic curve~\cite{FKL2}.

Several independent tests are applied to our examples,
and  confirm
our results in all cases.  This gives us confidence to claim that
our approximate functional equation method does indeed find
$L-$functions without the need to
first identify the underlying object.

Our data will be made available at
http://www.LMFDB.org/L/degree3  and 
/degree4.

\section{Prior methods}
We describe methods which have been used to find automorphic forms on $GL(2)$,
as well as recent approaches for higher rank groups.

\subsection{Holomorphic modular forms and Fourier series methods}
There are two types of modular forms on a Hecke congruence
group $\Gamma_0(N)\subset SL(2,\Z)$: holomorphic modular
forms and Maass forms.  Constructing holomorphic modular forms
is relatively straightforward.  If the level, $N$, is small, then one
can find an explicit basis of the ring of modular forms
on~$\Gamma_0(N)$.  For any level $N$ and any weight $k$, one can
use modular symbols~\cite{St1} to produce a
spanning set for the space $S_k(\Gamma_0(N))$.  Diagonalizing with
respect to the Hecke operators and the Atkin-Lehner operators
gives the basis of newforms, which have algebraic integer
Fourier coefficients.
This process has been completely automated~\cite{sage}.

For Maass forms, the situation is quite different.  Except for a thin
set arising from quadratic fields, there is no known
explicit construction of these functions.  Thus, for most
Maass forms one must rely on computer calculations to determine
numerical approximations~\cite{FL,Hej,Str}.  Recently,
such computations have been proven to be correct~\cite{BStV}.

The methods used to find Maass forms on $GL(2)$ make use of the Fourier
expansion and the transformation properties under elements of the group.
For $GL(d)$, $d\ge 3$, the Maass forms are functions of $\frac12(d^2+d-2)$
real variables.  Their Fourier expansions~\cite{Bu,Gol}
are $(d-1)$-fold sums over the integers, with an embedded sum over
cosets of $SL(d-1,\Z)$.  The complicated nature of the sum makes it
difficult to directly find Maass forms on $SL(3,\Z)$,
although partial success has been reported by Bian and Mezhericher~\cite{Bian, Mez}.
For $SL(4,\Z)$ and higher, that approach is probably infeasible.

\subsection{$L-$function methods}
A natural alternative to dealing with the Maass form is to work directly with
its $L-$function.  There is
no loss of information passing between those functions
since the spectral data and Fourier coefficients of the Maass form can be
recovered from the $L-$function. 
This approach was used successfully by Bian~\cite{Bian,Bian2},  
who made use of the 
functional equations for twists of the $L-$function, invoking the $GL(3)$
converse theorem.
The method involved solving a system of linear equations
with approximately 10,000 unknowns.
More details on Bian's method are given in Section~\ref{sec:other}.
Another approach was suggested by
Miller~\cite{Mil}, based on the Voronoi summation formula.

In this paper we present an approach which also involves finding the $L-$function, but we assume only one
functional equation and do not make the assumption of a functional
equation for twists.  Instead, we make use of the Euler product to
produce additional relations.  This approach enables us to work with a very
small system of equations, typically fewer than 30 unknowns,
although it is a non-linear system.

In the next section we describe the results we have obtained with our
method, and in Section~\ref{sec:thealgorithm} we give details about
our algorithm.

\section{Notation and Results}\label{sec:notation}
There have been several axiomatic definitions of ``$L-$function,''
most notably that of Selberg~\cite{Sel} who conjectured that
the functions in the Selberg class coincide with those arising
from automorphic representations.  Below, we give a more restrictive
set of axioms, which describes the properties that
are known or conjectured to hold for $L-$functions associated
to a unitary cuspidal automorphic representation of~$GL(d)$.

\subsection{$L-$function axioms}\label{sec:axioms}

We consider $L-$functions given by a Dirichlet series
\begin{equation}\label{eqn:DS}
L(s)=\sum_{n=1}^\infty \frac{a_n}{n^s},
\end{equation}
where $a_n\ll n^\epsilon$ for any $\epsilon>0$.  The estimate on the
Dirichlet coefficients is known as the \emph{Ramanujan bound}.

The Ramanujan bound implies that the Dirichlet series~\eqref{eqn:DS}
converges absolutely for $\sigma=\Re(s)>1$. We assume that $L(s)$ continues
to an entire function. (In general an $L$-function can have a pole at $s=1$,
but the $L$-functions we consider will be entire.)
We assume that the completed $L-$function, $\Lambda(s)$, is entire
and bounded in vertical strips, and
satisfies a functional equation of
the form
\begin{align}\label{eqn:FE}
\Lambda(s) =\mathstrut & N^{\frac{s}{2}} 
\prod_{j=1}^{d_1} \Gamma_\R(s+ \mu_j)
\prod_{j=1}^{d_2} \Gamma_\C(s+ \nu_j)
\cdot L(s)\cr
=\mathstrut & \varepsilon \overline{\Lambda}(1-s),
\end{align}
where $N$ is a positive integer called the \emph{level}, and $|\varepsilon|=1$.
Here $\Gamma_\R(s) = \pi^{-s/2} \Gamma(s/2)$
and $\Gamma_\C(s) = 2 (2 \pi)^{-s} \Gamma(s)$ where $\Gamma$ is the
Euler $\Gamma$-function, and $\overline{\Lambda}(z)=\overline{\Lambda(\overline{z})}$.
The integer $d = d_1 + 2 d_2$ is called the \emph{degree} of the $L-$function.
The analogue of the Selberg eigenvalue conjecture is that
$\Re(\mu_j),\Re(\nu_j)\ge 0$.

We also assume an Euler product of the form
\begin{equation}\label{eqn:EP}
L(s)= \prod_{p|N} L_p(s) \prod_{p\nmid N} L_p(s),
\end{equation}
where $L_p(s)=f_p(p^{-s})^{-1}$ with $f_p$ a polynomial
with $f_p(0)=1$.
If $p\nmid N$ then  we assume $f_p$ has degree~$d$ and satisfies
the self-reciprocal condition
$f_p(z) = \chi(p) (-1)^d z^d\, \overline{f_p}(z^{-1})$, where $\chi$
is a Dirichlet character mod~$N$, known as the
\emph{central character} of the $L$-function.
If $p|N$ then we assume $f_p$ has degree~$\le d-1$.

An equivalent formulation of the Ramanujan bound is that 
if $p\nmid N$ then all roots
of $f_p$ lie on the unit circle, and if $p|N$ then all roots of $f_p$ lie on
or outside the unit circle.  
This implies a more precise form of the Ramanujan bound:
\begin{equation}\label{eqn:preciseramanujan}
|a(p^m) | \le \genfrac{(}{)}{0pt}{}{m+d-1}{d-1} .
\end{equation}

\subsection{$L-$functions of Maass forms}\label{sec:SLNZ}
If $F$ is a Maass form on $SL(d,\Z)$ with Fourier coefficients
$(a_{n_1,\ldots,n_{d-1}})$, then its $L-$function $L(s,F)$
has Dirichlet coefficients
$a_n=a_{n,1,\ldots,1}$.  The $L-$function has an Euler product of the
form \eqref{eqn:EP} with $\omega_p=1$.
The Ramanujan bound has not been proven in this case.
The $L-$function satisfies a functional equation of the form
\eqref{eqn:FE} with $N=1$, $\varepsilon=1$, $d_1=d$, $d_2=0$,
and $\{\mu_1,\ldots,\mu_{d-1}\} =
\{ i\lambda_1,\ldots ,i\lambda_{d-1}\}$, where the $\lambda_i$
have simple expressions in terms of the eigenvalues
of the Maass form under the generators of the ring of invariant differential
operators, and $\mu_d=-\sum_{j<d} \mu_j$.  
The Selberg eigenvalue conjecture asserts that the $\lambda_j$ are real.
This was proven by Miller~\cite{Mil1} when $N=1$ and $d=3$.
See~\cite{Bu, Gol} for more details about these $L-$functions.

\subsection{Results for $SL(3,\Z)$}  
We continue the notation of Section~\ref{sec:SLNZ} with $d=3$.

By permuting the $\lambda_j$ and possibly replacing the $L-$function by
its dual (conjugate), we can assume $\lambda_1 \ge \lambda_2 \ge 0$
and $\lambda_3=-\lambda_1-\lambda_2$.
Thus, we specify the functional equation of the $L-$function,
equivalently the eigenvalues of the Maass form, by a point
$(\lambda_1, \lambda_2)$ 
below the diagonal of the first quadrant of~$\R^2$.  
The symmetric squares of
$SL(2,\Z)$ Maass forms have $\lambda_2=0$.

Weyl's law for $SL(3,\Z)$, proved by Miller~\cite{Mil1}, asserts that
\begin{equation}
\# \{(\lambda_1,\lambda_2) \ :
	\ \lambda_1^2 + \lambda_1\lambda_2 + \lambda_2^2 \le T\}
\sim \frac{vol(X)}{(4\pi)^{\frac52} \Gamma(\frac72)} T^{\frac52},
\end{equation}
where $X= SL(3,\Z)\backslash SL(3,\R)/SO(3,\R)$.
Miller also proved the Selberg eigenvalue conjecture for $SL(3,\Z)$
in the form $\lambda_1^2 + \lambda_1\lambda_2 + \lambda_2^2 \ge 80$.
(Miller proved slightly more, which is described in Figure~\ref{fig:SL3Z} below.)

The purpose of this paper is to describe a new method for 
finding $L-$functions satisfying
the axioms given in Section~\ref{sec:axioms}.
Our results fall into two categories: an extension of Miller's result that
certain regions do not contain the eigenvalues of any $SL(3,\Z)$ Maass forms,
and the experimental determination of a large number of $L-$functions
that appear to arise from $SL(3,\Z)$ 
Maass forms.
Both types of results are displayed in Figure~\ref{fig:SL3Z}.
\begin{figure}[htp]
\begin{center}
\scalebox{1.0}[1.0]{\includegraphics{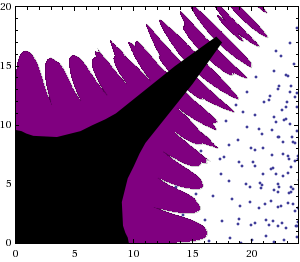}}
\caption{\sf 
The axes are $\lambda_1$ and $\lambda_2$ and the dots indicate
$(\lambda_1,\lambda_2)$ for the parameters in the
functional equation of $SL(3,\Z)$ $L-$functions,
for the first few examples we found.  The darkest region
is Miller's excluded region~\cite{Mil1}, and we prove,
assuming the Ramanujan conjecture, that the lighter shaded region contains
no eigenvalues.
} \label{fig:SL3Z}
\end{center}
\end{figure}

The dark-shaded region in Figure~\ref{fig:SL3Z} was shown by
Miller~\cite{Mil1}, using the explicit formula, to contain no
eigenvalues $(\lambda_1,\lambda_2)$ of any $SL(3,\Z)$ Maass forms.
We use the approximate functional
equation to prove a stronger result.

\begin{theorem}\label{thm:excluded} Assume the Ramanujan bound $|a_p|\le 3$.
There are no $L-$functions
with functional equation~\eqref{eqn:FE} with $d=3$ and $N=1$,
therefore
no $SL(3,\Z)$ Maass forms, with spectral parameters $(\lambda_1,\lambda_2)$
in the
lightly shaded region shown in Figure~\ref{fig:SL3Z}.
\end{theorem}

The proof, which is a new application of the approximate functional
equation, is given in Section~\ref{sec:excludedproof}.

Also shown in Figure~\ref{fig:SL3Z} are the 124 lowest
(in lexicographical order)
pairs $(\lambda_1,\lambda_2)$ of \hbox{$L-$functions} with functional
equation 
and Euler product as described in Section~\ref{sec:SLNZ} with $d=3$,
which were found in our search.
We believe that the figure
shows all cases with $\lambda_1\le 20$,
and may include all with $\lambda_1\le 24$.
Numerical values for the first 15 examples, and their first
two Dirichlet coefficients, are listed in Table~\ref{tab:SL3Z}.

Ce Bian (personal communication) has confirmed several dozen of our
examples and computed a few thousand Dirichlet coefficients in
each of those cases.

We have found more than 2000 spectral parameters 
$(\lambda_1,\lambda_2)$ and associated
Dirichlet coefficients.  These will be made available
at http://www.LMFDB.org/L/degree3.

\begin{table}
\[
\begin{array}{r|r||r|r}
\multicolumn{1}{c|}{\lambda_1} &
\multicolumn{1}{c||}{\lambda_2} &
\multicolumn{1}{c|}{a_2} &
\multicolumn{1}{c}{a_3 }
\\
\hline\hline
13.59658451 & 4.76468206 & 1.04846245 - 0.37523963 i & -0.49904094 - 0.27897508 i \\ 
14.14163558 & 2.38038848 & -0.10524097 + 0.75072694 i & 1.23599391 - 0.03911217 i \\ 
15.31863407 & 4.19173391 & 0.33541208 + 0.34590083 i & -0.27260607 - 0.97283444 i \\ 
 \hline 
15.74069912 & 7.85232504 & -0.39016397 - 0.62517821 i & 0.75666087 - 0.16479552 i \\ 
16.05436164 & 1.98365457 & -0.85221921 + 0.41800410 i & 0.67428665 - 0.33153749 i \\ 
16.40312474 & 0.17112189 & -0.42168648 + 1.06796797 i & -0.76802216 - 1.31329241 i \\ 
 \hline 
16.75957309 & 4.03941275 & 0.43484942 - 0.64408778 i & 0.67804925 - 0.27020173 i \\ 
17.34401833 & 7.13419625 & 0.55207716 - 0.70388045 i & 0.03833728 + 0.25697818 i \\ 
17.42523780 & 5.86543671 & -0.86407728 - 1.64923555 i & -0.22106187 - 0.90868355 i \\ 
 \hline 
17.50092882 & 1.90552951 & -0.28445287 + 0.01165958 i & 0.42424909 + 0.46219014 i \\ 
17.65047391 & 7.30146768 & -0.70569544 + 0.37307990 i & -0.29251139 + 1.07644911 i \\ 
17.86523805 & 10.34945597 & -0.86106525 + 0.35416339 i & -0.13894817 + 0.24569291 i \\ 
 \hline 
18.06174615 & 8.34053395 & -1.15143879 - 1.75641008 i & -0.49129679 + 0.21151351 i \\ 
18.08927092 & 5.30177649 & 0.32185117 + 0.32366378 i & -0.36219655 - 0.00596403 i \\ 
18.19481530 & 0.49254810 & -1.15442278 + 1.90981460 i & 0.33665482 - 0.10748733 i \\ 
 \hline 
\end{array}
\]

\caption{\sf \label{tab:SL3Z}
Spectral parameters $(\lambda_1,\lambda_2)$ and
sample Dirichlet series coefficients
for the first 15  $L-$functions of $SL(3,\Z)$ Maass forms. The numbers are truncated 
to 8 correct decimals. 
}
\end{table}

%
%
%

\subsection{Subgroups of $SL(3,\Z)$}\label{sec:SL3subgroups}

We applied the approximate functional equation method to
find $L-$functions associated to
Maass forms on subgroups of $SL(3,\Z)$.  In terms of the 
functional equation,
we have $d_1=3$, $d_2=0$, and $N>1$.
For $p|N$, the local factor in the 
Euler product \eqref{eqn:EP} will have degree at most~2.
In our calculations, we assume an arbitrary degree-2
local factor at the primes dividing the level:
$L_p(s) = (1+A_p p^{-s} + B_p p^{-2s})^{-1}$.  It may happen that
$B_p$ and/or $A_p$ is zero, if the local factor had degree 0 or~1.
By making no assumptions on the local factors, 
we have the opportunity for an independent check on
our calculations.

We applied the approximate functional equation method to
search for  $L-$functions of levels $N=2$ and $3$,  
and trivial central character,
and our search found nothing.
This seemed at first to be
problematic because when doing numerical experiments it can be
difficult to distinguish a true non-result from a non-result
due to a error in the computer program.  In this particular case,
it turns out that we experimentally discovered the theorem,
which
is relatively straightforward
but does not seem to be in the
literature: there are no Maass forms on $GL(3)$ which give
rise to $L-$functions of
squarefree level and trivial central character.
The most natural small index subgroups of
$SL(3,\Z)$ have outer automorphisms (analogous to the Fricke
involution for $\Gamma_0(N)$) which give level $p^2$ in the
functional equation.

So we searched for $L-$functions
of levels  $N=4$, and $9$ and found examples in each case.
For square level, the corresponding subgroups have outer automorphisms of 
order~3. Therefore, the signs which appear in the functional equation 
should be
3rd roots of~1, and we found examples of all three signs.
For level~$p^2$, in every case
the local factor we found at $p$ was
$L_p(s) = (1-\overline{\varepsilon} p^{-1-s})^{-1}$,
where $\varepsilon$ is the sign of the functional equation.
Note that we assumed a completely general degree 2 local factor,
but then found that it actually had degree~1.

Ralf Schmidt (personal communication) informed us that this
is the local factor from the Steinberg representation for~$GL(3)$,
which has level~$p^2$.  This is a strong independent
check of the validity of our calculations.

The first few $L-$functions for level~4  are given in 
Table~\ref{tab:level4}. 
The expected lifts from $GL(2)$
were found in our search.

\setlength{\extrarowheight}{4pt}
\begin{table}
\[
\begin{array}{c||r|r||c|r}
\multicolumn{1}{c||}{\epsilon} &
\multicolumn{1}{c|}{\lambda_1} &
\multicolumn{1}{c||}{\lambda_2} &
\multicolumn{1}{c|}{a_2 } &
\multicolumn{1}{c}{a_3}
\\
\hline\hline
1 & 8.23979752 & 2.64122672 & 1/2 + \delta & 1.06604225 + 0.48105885 i \\ 
1 & 10.64614640 & 5.40948610 & 1/2 + \delta & 0.06047305 - 0.43014710 i \\ 
1 & 11.72327411 & 2.02460527 & 1/2 + \delta & 0.59406181 - 0.27550311 i \\ 
1 & 12.42065536 & 4.72210152 & 1/2 + \delta & 0.99900525 - 0.54164255 i \\ 
 \hline 
 e(\frac13) & 9.63244452 & 1.37406028 & e(\frac23)/2 + \delta & 0.15012282 + 0.78657327 i \\ 
 e(\frac13) & 10.69551431 & 3.18301083 & e(\frac23)/2 + \delta & 0.62468273 - 0.09470423 i \\ 
 e(\frac13) & 11.25849528 & 1.13707825 & e(\frac23)/2 + \delta & -0.11242070 + 0.10083473 i \\ 
 e(\frac13) & 12.64900811 & 0.92183852 & e(\frac23)/2 + \delta & -0.80493907 - 0.12717437 i \\ 
 \hline 
e(\frac23) & 8.95466251 & 2.93659153 & e(\frac13)/2 + \delta & 0.36025530 - 0.35631998 i \\ 
e(\frac23) & 10.06466234 & 4.51670022 & e(\frac13)/2 + \delta & -0.90149321 - 0.57545547 i \\ 
e(\frac23) & 10.73782624 & 2.37415775 & e(\frac13)/2 + \delta & -0.53488262 + 0.88247158 i \\ 
e(\frac23) & 10.82653286 & 0.77122693 & e(\frac13)/2 + \delta & 0.58740526 - 0.33415990 i \\ 
\hline
\end{array}
\]

\caption{\sf \label{tab:level4}
Spectral parameters $(\lambda_1,\lambda_2)$, sign ($\epsilon$), and
sample Dirichlet series coefficients
for $L-$functions of level $4$ of degree $3$ Maass forms.
Here $e(x)=\exp(2\pi i x)$ and 
$\abs{\delta}<10^{-9}$. The numbers are truncated to what we believe are 8~correct
decimals.
}
\end{table}

\setlength{\extrarowheight}{0pt}

\subsection{Results for $Sp(4,\Z)$ and $GL(4,\Z)$}
We have also used our method to find degree~4 $L-$functions, that is,
$L-$functions with functional equation \eqref{eqn:FE}
and Euler product~\eqref{eqn:EP} with $d=4$.  At present, we
consider the case of level $N=1$; we will address the situation
of higher level in a subsequent paper.

\begin{figure}[htp]
\begin{center}
\scalebox{1}[1]{\includegraphics{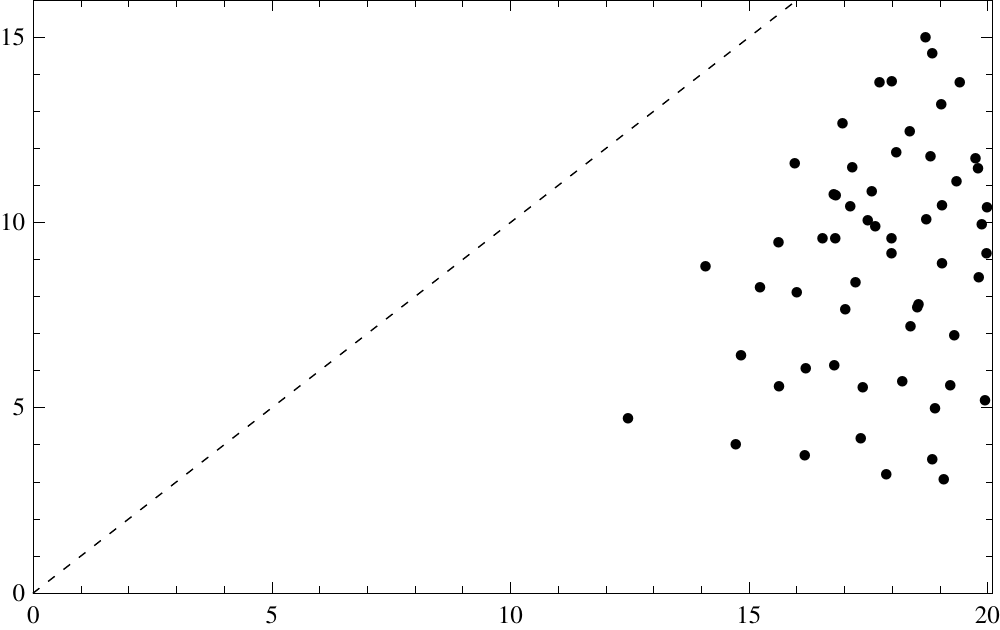}}

\caption{\sf
The axes are $\lambda_1$ and $\lambda_2$ and the dots indicate
$(\lambda_1,\lambda_2, -\lambda_1, -\lambda_2)$ 
for the first few $Sp(4,\Z)$ $L-$functions. 
Note that we may assume $\lambda_1\ge\lambda_2$.
} \label{fig:Sp4Z}
\end{center}
\end{figure}

It is natural to distinguish the special case 
where $\lambda_1=-\lambda_3$
and $\lambda_2=-\lambda_4$, and the Dirichlet coefficients $a_n$ are real.
We refer to these as ``$Sp(4,\Z)$ $L-$functions'' because we believe 
they are the spin $L-$functions of
the analogue of Maass forms on $Sp(4,\Z)$.  Since the functional equation
depends on two real parameters, we can provide a plot of our results, analogous
to the case of $SL(3,\Z)$.  Figure~\ref{fig:Sp4Z}
shows the first 58 examples that we have found for
$Sp(4,\Z)$ Maass form $L-$functions. The first 15 of these are also listed
in Table~\ref{tab:Sp4Z} with spectral parameters and the first four prime coefficients.

In Figure~\ref{fig:Sp4Z} one observes that the some of the first few $L-$functions
lie approximately along certain lines.  A similar phenomenon can be seen
in Figure~\ref{fig:SL3Z}. 

\begin{table}
\[
\begin{array}{r|r||r|r|r|r}
\multicolumn{1}{c|}{\lambda_1} &
\multicolumn{1}{c||}{\lambda_2} &
\multicolumn{1}{c|}{a_2} &
\multicolumn{1}{c|}{a_3} &
\multicolumn{1}{c|}{a_5} &
\multicolumn{1}{c}{a_7}
\\
\hline\hline
12.46875226 & 4.72095103 & 1.34260324 & -0.18745190 & -0.0016279 & 0.228229 \\ 
14.09372017 & 8.83036460 & -0.31644276 & 0.05294768 & 0.0261060 & 0.952439 \\ 
14.72547970 & 4.02669357 & 0.73680321 & 0.99830784 & 0.2677819 & -0.031085 \\ 
 \hline 
14.83756247 & 6.41243964 & 0.18715348 & -0.21546510 & 0.5492948 & 0.582942 \\ 
15.23336181 & 8.25653712 & 0.88602292 & 1.76959583 & -0.1217794 & -0.179071 \\ 
15.62467181 & 9.46217647 & -1.72583459 & 1.11618173 & -0.6444860 & -0.009432 \\ 
 \hline 
15.63453910 & 5.57458793 & 2.30065434 & 0.16113045 & -0.4928333 & -0.410839 \\ 
15.96343790 & 11.59606479 & 0.17708238 & -0.36379926 & -1.3728032 & 0.233579 \\ 
16.00467179 & 8.12341090 & -0.59748302 & -1.00280803 & -1.5990447 & 0.398796 \\ 
 \hline 
16.17300694 & 3.71016260 & 0.60584003 & 0.74786457 & -1.1437477 & 0.469017 \\ 
16.19395973 & 6.05889766 & 0.04742685 & 0.01653540 & 1.6223066 & -0.080626 \\ 
16.54574513 & 9.57822088 & 0.27348521 & 0.97263132 & 1.4040048 & -0.955692 \\ 
 \hline 
16.77928785 & 10.78021484 & -1.50813572 & -0.48797568 & 0.1694122 & -0.471615 \\ 
16.78915308 & 6.16376335 & -1.29328034 & -1.51751705 & -0.0818331 & -0.784244 \\ 
16.81031209 & 9.58414788 & 0.22614271 & -0.89126171 & 0.0562885 & -0.906897 \\ 
 \hline 
\end{array}
\]

\caption{\sf \label{tab:Sp4Z}
Spectral parameters $(\lambda_1,\lambda_2)$ and
sample Dirichlet series coefficients
for 15 $L-$functions of $Sp(4,\Z)$ Maass forms. The numbers are truncated to what we believe are 8 correct decimals.
}
\end{table}

The first $Sp(4,\Z)$ $L-$function, which has
$(\lambda_1,\lambda_2) = ( 12.46875,4.72095)$ has the curious property
that its lowest zero on the critical line is at $\frac12\pm 14.496i$,
which is higher than the lowest zero of the Riemann zeta-function.
This is discussed further by Bober,~et~al.~\cite{Bo}.
As a check on these calculations, we also compute the standard (degree~5)
$L-$function using the data from the spin (degree~4) $L-$function
(see section~\ref{sec:other}).

We have also found more than 200 $L-$functions of $SL(4,\Z)$ Maass forms
that are not $Sp(4,\Z)$ $L-$functions.
We list the eigenvalues and first
few coefficients of the smallest example.  We 
believe that all
the given digits are correct:
\begin{align*}
(\lambda_1,\lambda_2,\lambda_3) =\mathstrut&
(13.048202385582, 6.726210585782, -4.367763225255) \cr
a_2 =\mathstrut&-0.1912480688849 -0.587500805369\,i\cr
a_3 =\mathstrut&-0.652904479858 +0.96315195085\,i \cr
a_4 =\mathstrut&-0.5309393169 +0.22471678899\,i \cr
a_5 =\mathstrut&-0.3125111187 -0.0112785467\,i. \nonumber
\end{align*}

Numerical values of the spectral parameters and $a_2$ for the first 15 examples
we found
are listed in Table~\ref{tab:SL4Z}. They are normalized so that
$\lambda_1\geq \lambda_2\geq \lambda_3\geq \lambda_4$ and 
$\abs{\lambda_1}<\abs{\lambda_4}$ and in lexicographical order. We don't
claim that this list necessarily contains the first 15 \hbox{$L-$functions} with
this ordering.
The full list, including the Dirichlet coefficients, will be made available
at http://www.LMFDB.org/L/degree4.

\begin{table}

\[
\begin{array}{r|r|r||r}
\multicolumn{1}{c|}{\lambda_1} &
\multicolumn{1}{c|}{\lambda_2} &
\multicolumn{1}{c||}{\lambda_3} &
\multicolumn{1}{c}{a_2}
\\
\hline\hline
13.04820238 & 6.72621058 & -4.36776322 & -0.19124806 - 0.58750080 i \\ 
13.74207582 & 8.91090555 & -6.67116480 & -0.38362321 + 0.88277671 i \\ 
14.94493797 & 10.17319896 & -7.63795646 & 0.50916811 + 0.08060143 i \\ 
 \hline 
15.82354814 & 5.98442035 & -3.32054209 & 0.06707773 - 2.25215645 i \\ 
15.84333613 & 8.14020386 & -6.39808715 & -0.93727868 + 0.45045257 i \\ 
15.94944916 & 11.22357307 & -7.46849155 & -0.01250803 - 0.16785927 i \\ 
 \hline 
15.99896071 & 5.55250290 & -4.28979199 & 0.76641492 + 0.72684214 i \\ 
16.04868870 & 7.44354895 & -5.74749580 & 0.39710090 - 1.01290587 i \\ 
16.17564264 & 9.31222108 & -7.99494686 & -0.25094605 - 1.40588826 i \\ 
 \hline 
16.72584930 & 10.57549478 & -8.05016487 & -1.81300469 - 0.42294306 i \\ 
16.89972715 & 2.27258771 & -6.03583588 & 0.55659019 + 0.92818484 i \\ 
16.92627200 & 12.47850050 & -10.63475738 & 0.33041049 + 0.59891916 i \\ 
 \hline 
16.95696382 & 9.53674609 & -7.32770475 & -0.66772182 + 0.75485571 i \\ 
17.16046064 & 5.51936319 & -3.67883841 & -0.31914164 - 0.69706013 i \\ 
17.27511800 & 7.18638875 & -4.38401627 & -0.13829583 + 1.09540674 i \\ 
 \hline 
\end{array}
\]

\caption{\sf \label{tab:SL4Z}
Spectral parameters $(\lambda_1,\lambda_2,\lambda_3)$
and $a_2$
for $L-$functions of $SL(4,\Z)$ Maass forms. The numbers are truncated to what we believe are 8 correct decimals.
}
\end{table}

\subsection{Other types of degree 4 $L-$functions}
As a test of the approximate functional equation method, we now apply it
to a different type of degree 4 $L-$function.
In this example, we search for $L-$functions
satisfying a functional equation of the form
\begin{equation}\label{eqn:FEsiegel}
\Lambda(s)=\Gamma_\C(s+\tfrac12) \Gamma_\C(s+k-\tfrac32) L(s) = \Lambda(1-s).
\end{equation}
That is the functional equation of the spin $L-$function of a Siegel modular
form of weight~$k$ on $Sp(4,\Z)$.  As a demonstration of the method,
let $k=20$.  It is known that there is exactly one such $L-$function,
associated to the cusp form denoted $\Upsilon_{20}$.  Its 
first few Dirichlet coefficients,
in the arithmetic normalization, are the following rational integers~\cite{Sko}:
\begin{align}
2^{37/2} a_2 =\mathstrut & -840960 \cr
3^{37/2} a_3 =\mathstrut & \phantom{\mathstrut - \mathstrut}346935960 \cr
4^{37/2} a_4 =\mathstrut & \phantom{\mathstrut - \mathstrut}316975677440 \cr
5^{37/2} a_5 =\mathstrut & -5232247240500 .
\end{align}


We applied the method of Section~\ref{sec:thealgorithm} to search for the Dirichlet coefficients
of an $L-$function satisfying the functional equation \eqref{eqn:FEsiegel} with $k=20$,
 and having
a degree-4 Euler product.
We used the first 90 terms of the Dirichlet series, which is 30 unknowns because
a degree-4 local factor is determined by $a_p$ and $a_{p^2}$.
We chose 30 triples $(\frac12+i t,g_1, g_2)$ with $0 \le t \le 33$, selected so that
the truncation error in eliminating $a_n$ for $n> 90 $ is  less than $10^{-12}$.
After using the Euler relations to convert to a non-linear system, we solved
the system 1000 times using the secant method with random starting values.
Only one solution was found (and it was found repeatedly), which had:
\begin{align}
2^{37/2} a_2 =\mathstrut &\mathstrut -840960.000014382993 \cr
3^{37/2} a_3 =\mathstrut &\mathstrut \phantom{\mathstrut - \mathstrut}346935960.202402521 \cr
4^{37/2} a_4 =\mathstrut &\mathstrut \phantom{\mathstrut - \mathstrut}316975677305.428055 \cr
5^{37/2} a_5 =\mathstrut &\mathstrut -5232247227434.34870. 
\end{align}
We see that the (rescaled) $a_2$ is convincingly close to an integer, and it is
the correct integer.  Replacing $a_2$
by its apparent exact value, and the re-solving for the other coefficients, we find:
\begin{align}
3^{37/2} a_3 =\mathstrut & \phantom{\mathstrut - \mathstrut}346935959.999827467 \cr
4^{37/2} a_4 =\mathstrut & \phantom{\mathstrut - \mathstrut} 316975677440.106735 \cr
5^{37/2} a_5 =\mathstrut & -5232247240497.56648 . 
\end{align}
Now the rescaled $a_3$ is convincingly close to an integer, and again it is the correct integer.
Replacing $a_3$
by its apparent exact value, and solving again we find:
\begin{align}
4^{37/2} a_4 =\mathstrut & \phantom{\mathstrut - \mathstrut}316975677439.993180 \cr
5^{37/2} a_5 =\mathstrut & -5232247240498.27543. 
\end{align}
This time, the situation is less clear.
The rescaled $a_4$ certainly appears to be close to an integer;
in fact the correct integer.
But we had truncated the Dirichlet series to give
12~digits of accuracy.  Since the true value of $4^{37/2} a_4 $
is an integer with 12 digits, our calculation should not have been adequate
to determine that integer exactly.  It is somewhat surprising
that we obtain a few more correct digits than we should have expected.

This example illustrates that the approximate functional equation method
should be applicable to other types of $L-$functions.  
We will address these extensions of the method in a subsequent paper.

\section{The method:  approximate functional equations}\label{sec:thealgorithm}

Our method of locating $L-$functions is based on
a smoothed approximate functional equation.
The approximate functional equation is a tool that typically is used
to prove results about $L-$functions (such as moments), or to
calculate $L-$functions~\cite{Rub}.  Here we turn it into a method of
locating $L-$functions by exploiting the fact that the
approximate functional equation contains a free parameter.
This approach was inspired by Michael Rubinstein's \emph{Lcalc}
program~\cite{Rub, Rub2}.

\subsection{Smoothed approximate functional equation}
Let
\begin{equation}
   L(s) = \sum_{n=1}^{\infty} \frac{a_n}{n^s}
\end{equation}
be a Dirichlet series that converges absolutely in a half plane, $\Re(s) > \sigma_1$.

Let
\begin{equation}
    \label{eq:lambda}
    \Lambda(s) = Q^s
                 \left( \prod_{j=1}^a \Gamma(\kappa_j s + \lambda_j) \right)
                 L(s),
\end{equation}
with $Q,\kappa_j \in {\mathbb{R}}^+$, $\Re(\lambda_j) \geq 0$.
Note that the functional equation  we present here is somewhat more
general that given in the previous section.
We assume that:
\begin{enumerate}
    \item  $\Lambda(s)$ has a meromorphic continuation to all of ${\mathbb{C}}$ with
           simple poles at $s_1,\ldots, s_\ell$ and corresponding
           residues $r_1,\ldots, r_\ell$.
    \item $\Lambda(s) = \varepsilon \cj{\Lambda(1-\cj{s})}$ for some
          $\varepsilon \in {\mathbb{C}}$, $|\varepsilon|=1$.
    \item For any $\sigma_2 \leq \sigma_3$, $L(\sigma +i t) = O(\exp{t^A})$ for some $A>0$,
          as $\abs{t} \to \infty$, $\sigma_2 \leq \sigma \leq \sigma_3$, with $A$ and the constant in
          the `Oh' notation depending on $\sigma_2$ and $\sigma_3$. \label{page:condition 3}
\end{enumerate}

To obtain a smoothed approximate functional equation with desirable
properties, Rubinstein \cite{Rub} introduces an auxiliary function.
Let $g: \C \to \C$ be an entire function that, for fixed $s$, satisfies
\begin{equation*}
    \abs{\Lambda(z+s) g(z+s) z^{-1}} \to 0
\end{equation*}
as $\abs{\Im(z)} \to \infty$, in vertical strips,
$-x_0 \leq \Re (z) \leq x_0$. The smoothed approximate functional
equation has the following form.
\begin{theorem}
    \label{thm:formula}
    For $s \notin \cbr{s_1,\ldots, s_\ell}$, and $L(s)$, $g(s)$ as above,
    \begin{align}
         \label{eq:formula}
         \Lambda(s) g(s) =
         \sum_{k=1}^{\ell} \frac{r_k g(s_k)}{s-s_k}
         + Q^s 
           \sum_{n=1}^{\infty} \frac{a_n}{n^s} f_1(s,n) 
         + \varepsilon Q^{1-s} 
                \sum_{n=1}^{\infty} \frac{\cj{a_n}}{n^{1-s}} f_2(1-s,n)
    \end{align}
    where
    \begin{align}
        \label{eq:mellin}
        f_1(s,n) &:= \frac{1}{2\pi i}
                    \int_{\nu - i \infty}^{\nu + i \infty}
                    \prod_{j=1}^a \Gamma(\kappa_j (z+s) + \lambda_j)
                    z^{-1}
                    g(s+z)
                    (Q/n)^z
                    dz \notag \\
        f_2(1-s,n) &:= \frac{1}{2\pi i}
                    \int_{\nu - i \infty}^{\nu + i \infty}
                    \prod_{j=1}^a \Gamma(\kappa_j (z+1-s) + \cj{\lambda_j})
                    z^{-1}
                    g(s-z)
                    (Q/n)^z
                    dz
    \end{align}
    with $\nu > \max \cbr{0,-\Re(\lambda_1/\kappa_1+s),\ldots,-\Re(\lambda_a/\kappa_a+s)}$.
\end{theorem}
Note that, despite its name, the approximate functional equation
is an exact formula for the $L-$function throughout the complex plane.

We assume $L(s)$ continues to an entire function, so the first
sum in \eqref{eq:formula} does not appear.  For fixed
$Q,\kappa,\lambda,\varepsilon$, sequence~$\{a_n\}$, and function $g(s)$,
the right side of \eqref{eq:formula}
can be evaluated to high precision. In the next section we describe
how we use this observation to model the functional equation.

\subsection{Modeling the functional equation}\label{ssec:checkFE}
The Z-function (or Hardy Z-function) of an \hbox{$L-$function} is an analytic
function which is smooth and real on the $\frac12$-line and satisfies
\hbox{$|Z(\frac12+it)| = |L(\frac12+ it)|$} for $t\in \R$.  In terms of
the completed $L-$function~\eqref{eq:lambda},
\begin{equation}\label{eqn:Z}
Z(\tfrac12+it) = (\varepsilon Q)^{-\frac12}
\prod_{j=1}^a |\Gamma(\kappa_j(\tfrac12+it)+\lambda_j)|^{-1} 
\cdot \Lambda(\tfrac12+it),
\end{equation}
when $t\in \R$.  That formula is not a valid definition of
$Z(\frac12+it)$ if $t\not\in \R$, but it is
sufficient for our purposes because we will only evaluate the $L-$function
on the critical line.  We will find the Z-function to be convenient
because all calculations will involve real numbers instead
of complex numbers.

Fix all the parameters in the functional equation of an $L-$function and
write \eqref{eq:formula} as
$\Lambda(s) g(s) = \AppFE(s,Q,\kappa,\lambda,\varepsilon,\{a_n\},g)$,
where we suppress the contribution of the poles because we will assume
our $L-$functions are entire.
Use~\eqref{eqn:Z} to rewrite~\eqref{eq:formula} as
\begin{align}\label{eq:rearranged}
Z(s) = Z(s,g) =\mathstrut &
g(s)^{-1}  (\varepsilon Q)^{\frac12}
\prod_{j=1}^a |\Gamma(\kappa_j s +\lambda_j)|^{-1} \,\AppFE(s,Q,\kappa,\lambda,\varepsilon,\{a_n\},g) \cr
=\mathstrut & g(s)^{-1}  \,\RHS(s,Q,\kappa,\lambda,\varepsilon,\{a_n\},g),
\end{align}
say.

The only free parameters in the above equation are the complex number~$s$
and the function~$g$.
The left side of \eqref{eq:rearranged} does not depend on the function~$g$.
Thus, for each fixed~$s$, if we evaluate the right side of \eqref{eq:rearranged}
with two different $g$ functions $g_1$, $g_2$, then setting them equal gives an
equation in the Dirichlet coefficients of the $L-$function.

Thus, if we fix $s_0=\frac12+i t_0$ and choose two functions
$g_1, g_2$ then
\begin{equation} \label{eq:linear}
 g_1(s_0)^{-1}  \,\RHS(s_0,Q,\kappa,\lambda,\varepsilon,\{a_n\},g_1)
- 
 g_2(s_0)^{-1}  \,\RHS(s_0,Q,\kappa,\lambda,\varepsilon,\{a_n\},g_2) = 0
\end{equation}
is a linear equation which the Dirichlet coefficients of the $L-$function must satisfy.

We can produce many equations by choosing various $(s,g_1,g_2)$ combinations.
If we have chosen functional equation parameters which correspond to a
genuine $L-$function, then its Dirichlet coefficients will be a solution to
that system.  If we have chosen parameters for which there is no $L-$function,
then it may be possible to detect this through an inconsistency in the system.

We expand on this approach in the next section.

\subsection{The method}\label{sec:method}
It is easiest to first think of the functional equation as given,
and focus on determining the Dirichlet coefficients.  For the case of
Maass forms, we do not actually know the functional equation in advance:
we must find the spectral parameters $\lambda_1,\ldots,\lambda_{d-1}$
which determine the functional equation.  This is an added complication
which we address in Section~\ref{sec:searching}.

Assume we are given the functional equation of an $L-$function
and we want to
find its Dirichlet coefficients.
The main idea is to first choose various
$s_j=\tfrac12+i t_j$, and pairs of test functions~$g_{j,1}$ and
$g_{j,2}$ to form the system of linear equations
$\{Z(s_j,g_{j,1}) - Z(s_j,g_{j,2}) = 0\}$.

\step{1}{Forming the linear system}

Typically we choose $g(s) = e^{-i \beta s + \alpha s^2}$ with $\beta,\alpha\in \R$.
This is an allowable
test function in the approximate functional equation provided
$\alpha>0$, or $\alpha=0$ and $|\beta|<\frac{\pi}{2}\sum \kappa_j$.

We illustrate the ideas with a numerical example.  The example
concerns the $L-$function of a Maass form on $SL(3,\Z)$, so in \eqref{eq:lambda} we
set $Q=\pi^{-3/2}$, $a=3$, $\kappa_1=\kappa_2=\kappa_3=\frac12$,
 and $\varepsilon=1$.  
This specifies all the parameters in the functional equation
except $\lambda_j$.
We illustrate below with the test case
$\lambda = (i \lambda_1,i \lambda_2, i \lambda_3) = (8.4 i,14.2 i,-22.6 i)$.
Those parameters were chosen randomly, so there is no reason to
expect that this corresponds to an actual $L$-function.
In fact, in Section~\ref{sec:excludedproof} we will prove that
there is no $L$-function with that functional equation.

Note that we must treat the real and imaginary parts of
the Dirichlet coefficients as separate real variables,
so we write
$a_n = b_1(n) + i b_2(n)$ with $b_1(n), b_2(n) \in \R$.

We apply the procedure described in the previous section.
This involves choosing a complex number~$s_0$ and two
test functions~$g_1$,~$g_2$.
Choose $s_0=\frac12+i$, and let $g_1(s)= e^{-i s/4}$.
  A numerical calculation finds,
\begin{align} \label{eqn:beta1}
g_1(s_0)^{-1}  &\,\RHS(s_0,Q,\kappa,\lambda,\varepsilon,\{a_n\},g_1)= \\
&\mathstrut 1.2860+7.7907\, b_1(2)+1.9245\, b_2(2) 
+0.6653\, b_1(3)-1.0163\, b_2(3) \cr
&+0.03319\, b_1(4) -  0.1197\, b_2(4) 
+\cdots  +4.17\cdot 10^{-6}\, b_1(8) -2.37 \cdot 10^{-7}\, b_2(8) + \cdots.\nonumber
\end{align}
The terms decrease rapidly: the coefficients of $b_{1}(15)$
and $b_2(15)$ are less than $2\cdot 10^{-13}$.  The exponential decay in the
contributions of the Dirichlet coefficients is a general
phenomenon~\cite{Book, FRS}.

With the same $s_0$, now choose $g_2(s)= e^{-i s/2}$.  This gives
\begin{align} \label{eqn:beta2}
g_2(s_0)^{-1}  &\,\RHS(s_0,Q,\kappa,\lambda,\varepsilon,\{a_n\},g_2)= \\
&\mathstrut 2.02524 +2.2117\, b_1(2)+1.8524\, b_2(2)
+0.7408\, b_1(3)+0.2210\, b_2(3) \cr
&+0.1019\, b_1(4) +0.0922\, b_2(4) 
+\cdots  2.66 \cdot 10^{-6}\, b_1(8) -0.0000143 \, b_2(8) + \cdots.\nonumber
\end{align}
Subtracting \eqref{eqn:beta1} from \eqref{eqn:beta2}, we obtain the
following equation which the Dirichlet coefficients of the
$L-$function must satisfy, \emph{if} such an $L-$function exists:
\begin{align} \label{eqn:beta12}
0=\mathstrut & 0.7392 - 5.5790 \, b_1(2) - 0.07216 \, b_2(2) 
+0.07541 \, b_1(3) +1.23745  \, b_2(3) \\
&0.06876 \, b_1(4)  +0.2119 \, b_2(4)  
+\cdots   - 1.51 \cdot 10^{-6}\, b_1(8) -0.0000141  \, b_2(8) + \cdots.\nonumber
\end{align}

We have just found \emph{one} equation which the Dirichlet coefficients
of the $L$-function (if it exists) must satisfy.
The next steps involve 
truncating the system so that we
are dealing with finitely many unknowns,
and choosing various $(s_0,g_1,g_2)$ combinations 
to form a system of equations.

\step{2}{Truncating the system}

The linear equation from Step 1 has infinitely many unknowns.
However, the coefficients of those unknowns decrease rapidly.
Thus, once we select a working precision, we can truncate the system
with an error less than our working precision.

For the remaining steps, we assume that a working precision, $\delta$,
has been selected (for example, $\delta=10^{-12}$),
and a given number of Dirichlet
coefficients, $J$, have been selected (for example, $J=30$, so we
need only consider $a_1,\ldots,a_{30}$).

\step{3}{Choosing the parameters}

Having selected a working precision $\delta$ and a number of
Dirichlet coefficients $J$,
there are limits to which $(s,g_1,g_2)$ combinations can be used.
Namely, we require that the error in the approximate functional equation
be less than $\delta$ when we use these particular values of $s$ and $g_j$ and omit those
terms $a_n$ with $n>J$.  Note that this requires the Ramanujan bound
or some other estimate on the Dirichlet coefficients.

In particular, given $\delta$ and $J$, for each function $g$, there
is a limited range of $s$  (which may be empty)
which meets the requirement of the chosen truncation error.  

In practice, often it is effective to fix $g_1$ and $g_2$,
and to vary $s_j = \frac12+i t_j$ to obtain different linear equations.
However, in order to keep the system of equations as well conditioned
as possible, one should aim to have the same rate of convergence in all
equations. This can be achieved by varying the functions $g_1$ and $g_2$
with the $s_j$.

\step{4}{Using the Euler product}

Except in some low-degree cases, the linear system will have
infinitely many solutions.  For example, if there are two $L-$functions
with the same functional equation, then clearly the linear system
will have infinitely many solutions.

In order to get around this, we make use
of the Euler product and convert to a non-linear system.
The difficulty of working with a non-linear system is partially
offset by the fact that the number of unknowns will be much smaller.
The Euler product immediately gives that all Dirichlet
coefficients are determined by the coefficients at prime powers.
Furthermore, for a degree-$d$ Euler product, the local factors
are determined by $[d/2]$ parameters.
For example, for degrees
2 and 3, all coefficients are determined by $a_p$.
And for degrees 4 and 5, all coefficients are determined by $a_p$
and $a_{p^2}$.

We apply the multiplicative relations $a_6=a_2 a_3$, $a_{12} = a_3 a_4$,
etc., to reduce the number of unknowns.  Then we apply the recursions
arising from the shape of the  local factors (for example,
$a_{p^2} = a_p^2-\overline{a_p}$ in the case of $SL(3,\Z)$
Maass forms) to obtain a system of equations with a minimal
number of unknowns. So, for example, for $SL(3,\Z)$
Maass forms the unknowns will be $\{a_p : p \text{ prime, } p\leq J\}$.

Note:  we refer to the $a_n$ as the ``unknowns,'' but in practice
we must treat the real and imaginary parts of $a_n$, denoted
$b_1(n)$ and $b_2(n)$, as separate real unknowns.  Since
$a_1=1$, the original linear system has $2(J-1)$ real unknowns.

\step{5}{Solving the non-linear system}

The non-linear system will have many solutions. Most of these
are of no interest and are simply accidental solutions.

We use the secant method to solve the system. This is done repeatedly,
typically between 100 and 1000 times,
with different starting points. The starting
points are always chosen to have values within the Ramanujan bound since
this is supposedly true for the interesting solutions. Typically 
0--5 different solutions are found. In most cases, the size of the $a_n$ in the
solutions will grow rapidly with $n$. For the accidental solutions, 
the growth is (almost) as fast as the rate of decrease in the coefficients.
For the solutions that are close to approximations
of true $L-$functions, this is not true. Here the growth in the size of the $a_n$ 
is smaller. 
Thus, the Ramanujan bound is an important ingredient to the method because
it allows us to eliminate unwanted solutions.
However, 
there are more efficient ways to avoid extraneous solutions, see 
Section~\ref{sec:searching}.

In practice, it turns out that it is better
to begin with not too many unknowns in order to find 
any of the solutions at all. Once a solution has been found, the number of
unknowns may be increased and a new larger system is solved with the
previous  solution as a starting point.
Typically we start with $10\leq J\leq40$ and
end with $60\leq J\leq150$.

\medskip

The above procedure describes a method for finding the Dirichlet coefficients
of an \hbox{$L-$function}, given its functional equation and the shape of its Euler product.
In the case of Maass form $L-$functions, there is the added difficulty of finding
the spectral parameters which determine the functional equation.  This is discussed
in the next section.

\section{Searching for the spectral parameters}\label{sec:implementation}

When searching for $L-$functions of Maass forms
on $SL(d,\Z)$ or one of its subgroups, the spectral parameters
in the functional equation are not known in advance.
The greatest
challenge is to find good approximations to these parameters.

\subsection{Initial scan}\label{sec:searching}

The first step is to fix a box in which to initially search for spectral parameters
$\lambda_j$ corresponding to an $L-$function. Typically, we use boxes of
side length $2$, so for $SL(3,\Z)$ and $Sp(4,\Z)$ this will be a
square (with side $2$), and for $SL(4,\Z)$, a cube. In this box
we make a grid of equally spaced points with a typical step size
of $1/5$ to $1/25$. 

Next, we choose the number of Dirichlet coefficients, $J$, to use
and set the working precision (estimated maximal size of truncation error)
to $\delta=0.5\cdot 10^{-4}$. We fix test functions $g_1$ and $g_2$ of the
form $g_i(s)=e^{-i\beta_i s+\alpha s^2}$. After these choices are made,
we determine the range $t\in[T_1,T_2]$ of $s=\frac12+i t$ on  the critical line for
which we have a truncation error less than $\delta$ in the 
approximate functional equation. If this range is too small, we increase
$J$ until we achieve a large enough range. Typically the range $[T_1,T_2]$ will
be approximately  $[2\min\{\lambda_j\},2\max\{\lambda_j\}]$. 

Now we determine the number of unknowns, $m$, of the 
non-linear systems of 
equations. This will depend on whether we assume the $L-$functions are self-dual
(meaning that the Dirichlet coefficients are real), and on the degree of the
local factors of the Euler product. In the case of $SL(3,\Z)$, $m$ 
is twice 
the number of primes less than or equal to~$J$. 
For $Sp(4,\Z)$, it is the number of primes or squares of primes less than or equal to~$J$. 
We choose
$m+k$ (typically with $k=8$) values of 
$s_i=\frac12+t_i$ in the range 
$T_1 < t_j < T_2$. Of these points, $m$ of them will be used to make
the systems of $m$ non-linear equations in $m$ unknowns, and the other
$k$ will be used to detect candidates for first approximations to the spectral 
parameters, as described in the next subsection.

Now for each point in our grid of spectral parameters $\lambda_j$,
we go through the following procedure.
Since the parameters $\lambda_j$ and test functions 
$g_1$ and $g_2$ are fixed, each of the
points $s_i$ gives a linear equation of the form~(\ref{eq:linear}). 
The Euler product is used to turn these into $m+k$ non-linear
equations with $m$ unknowns. We choose $m$ of these equations (same
for all points in the grid) and then use the secant method to solve
this non-linear system of equations. The system of equations is
solved a fixed number of times (typically $50$) with starting points
chosen randomly
each time.
In addition to these random starting points, we also use solutions
that have already been found at adjacent points in the grid as 
additional starting values.
Sometimes the system with a given starting point will converge to a 
solution with given accuracy
within a set number of steps and sometimes it will not. For some points in the
grid it will converge every time and for some never. These solutions
are considered ``the same'' if the first few coefficients agree to a certain
accuracy. So finally we will have a number (typically 0--5) of different
solutions at each points in the grid.

\subsection{Refining the search}
Once we have a number of different solutions for each point in the grid,
we can find candidates for the functional equation parameters.
This is when the $k$ extra equations come into play. 
Suppose the spectral
parameters used in the system of equations are close to those of an $L-$function
and the Dirichlet coefficients of the solution are also close
to those of the $L-$function. Then those coefficients should also
be close to solutions to the $k$ extra equations.
Hence we want to find spectral parameters for which there
is a solution that makes the $k$ extra equations close to consistent.
We do this by substituting the solution found at each point into the extra equations,
forming ``indicators'' which should measure how far the point is from a true $L-$function.
We then interpolate to determine where the indicators approximately vanish.

\subsubsection{Collecting nearby solutions}
Since there may be several solutions at each point, some of which may
be extraneous, we first collect those solutions which are likely
to correspond to the same $L-$function.

Fix a point $P$ in the grid 
and choose $D$ adjacent points in the grid where $D$ is the dimension of the
grid. In the case of $SL(3,\Z)$, we have $D=2$ and choose the points just
above and just to the right of $P$. For each of the solutions at $P$ we 
identify the solutions at the adjacent points that are closest. Here
`closest' in the
sense that the first few coefficients differ by as little as possible.
Now we have sets of related solutions at $D+1$ adjacent points,
and each set could potentially correspond to an $L-$function. 

\subsubsection{Comparing nearby solutions}
Now we use each set of solutions and the $k$ extra equations to form indicators, and then interpolate to
determine where those indicators (approximately) vanish.

Choose $D$ of the $k$ extra equations at each of the
$D+1$ points. Substitute the related solutions of a set into the $D$ extra 
equations to get $D$ complex numbers (indicators) at each of the $D+1$ points.
Using linear interpolation we find a point where the  indicators are approximately~$0$,
which gives us a candidate for an $L-$function.
This is done for each of the
$\binom{k}{D}$
choices of extra equations, so
we get 
$\binom{k}{D}$
candidates. If a large enough proportion (50\%) of these
candidates is close enough to $P$ this is considered a good candidate.
We evaluate the average (of the central 50\%) of these candidates and consider
this to be an initial candidate for the spectral parameters of an $L-$function,
which is passed on the the next step of the process.
We also save an interpolation of the 
first few Dirichlet
coefficients of the solutions in the set. These coefficients will serve 
as a starting point when trying to improve the precision of this candidate.

This process is repeated for each of the related sets of solutions at $P$ and
then for each of the points in the grid. The result of this process is
a list of initial candidates for spectral parameters of $L-$functions
with a first approximation of its Dirichlet coefficients. 

\subsection{Zooming in on a candidate}
The final step is to try, for each of the initial candidates, to zoom in to a
good approximation for the spectral parameters. This is done in steps where
the number of coefficients is gradually increased and hence also the size
of the truncation error decreased. Each of these steps is similar to the 
linear interpolation that gave the initial candidate. A system is created and solved
at $D+1$ adjacent points with the approximation of the Dirichlet coefficients
as a starting point. If there are solutions at each of the points, then we get
$\binom{k}{D}$
new approximations in the same way as before. If there is no solution
or if the
average of the new approximations is far from the first then we 
stop pursuing that possible $L-$function, otherwise we continue.
This is repeated with increasingly smaller steps and more
coefficients. This process is deemed a success when the maximum difference among the
$\binom{k}{D}$
approximations is less than a given bound (we  typically require $8$ correct
decimals). The average (of the central 50\%) of these approximations is saved
as a highly probable candidate for spectral parameters of an $L-$function.

\medskip

In practice, this process converges very quickly, typically yielding one
or two additional decimal digits in each step.
We have implemented the above procedure in Mathematica and will make our data available
online~\cite{lmfdb}.

\section{Checks on the results}\label{sec:other}

In this section, we discuss some independent checks which support
the claim that our procedure
finds approximations to $L-$functions of Maass forms.
 
\subsection{Functional equation for twists}
Bian~\cite{Bian2} computed the first four examples of \hbox{$L-$functions} of
Maass forms for $SL(3,\Z)$. 
His method uses a similar starting point, but he also
uses the functional equation of many twists of the $L-$function. This has the
drawback of requiring a large number of coefficients. Already for the first few
examples, his method requires a linear system of equations with up to
$10000$ unknowns.
Our search found his four examples and a comparison of spectral parameters
and Dirichlet
coefficients finds agreement to within the claimed precision. Our 
method gives a higher precision
in the spectral parameters and the first few coefficients, but 
Bian's method produces
more coefficients. 
Finally, Bian has verified the first $50$ $SL(3,\Z)$ \hbox{$L-$functions}
found in our search.

\subsection{Properties of the coefficients}
Another check is that in our method we leave $a_{32}$ (and also $a_{64}$) as unknowns.
Hence we obtain a computed value of $a_{32}$ that we can compare with the 
the value derived from the computed value of $a_2$. For all
our examples for $SL(3,\Z)$ these values agree to at least three
decimal places 
and in some cases (with larger
spectral parameter requiring more coefficients) to as many as eight decimal places. This
also gives an independent estimate for the accuracy of the coefficients,
which compares well with our
estimate based on differences in the approximations in the computation.

A strong check involving the coefficients is that
the local factors at the primes dividing the
level for $GL(3)$ $L-$functions agree to high precision with the
predicted value.
This is discussed in Section~\ref{sec:SL3subgroups}.


\subsection{Lifts}
A check which we consider to be particularly
strong is the analytic continuation and functional equation for higher
degree lifts of our $L-$functions.

For various $GL(3)$ 
examples,
we computed the symmetric square of the $L-$functions,
which has degree~6.  
That is, we used the functional equation and Dirichlet coefficients
of the degree~3 $L-$function to determine the functional equation and
Dirichlet coefficients of the degree~6 $L-$function.
We then used the method in Theorem~4.2 of~\cite{FRS}
to verify that the $L-$function satisfies the conjectured functional
equation.  In all cases we find that degree~6 $L-$function satisfies
the expected functional equation.

For the $GL(4)$ $L-$functions, the symmetric square has degree~10.
Such high degree \hbox{$L-$functions} require a large number of
coefficients to be computed accurately, so we were not able to
obtain convincing results using the symmetric square.
But for the specific case
of $Sp(4,\Z)$ there is another option.  We claim that 
we have found the spin $L-$functions associated to Maass forms on $Sp(4,\Z)$.
If that is true, then we can use the functional equation data and
Dirichlet coefficients to compute the standard $L-$function, which
has degree~5.
(That
the standard $L-$function can be reconstructed from the
spin $L-$function can be seen from an examination
of the explicit representations given in~\cite{FRS}.)
We computed the degree~5 \hbox{$L-$functions} 
in several cases and found that they
satisfy the expected functional equation to high accuracy.
We consider this to be a strong confirmation of the
legitimacy of our calculations.  Numerical data will be made
available at http://www.LMFDB.org/L/degree5.

\section{Proof of Theorem~\ref{thm:excluded}}\label{sec:excludedproof}
In this section we describe the proof that the lighter shaded region in Figure~\ref{fig:SL3Z}
contains no degree 3 $L-$functions.  The idea of the proof is to use the 
approximate functional equation method
of Section~\ref{sec:thealgorithm}  to form a system of linear equations in the
Dirichlet coefficients, and then verify that the solution to that system
violates the Ramanujan bound, when $(\lambda_1,\lambda_2)$ are in the shaded region.

We will write out the details for the point
$(\lambda_1,\lambda_2) = (8.4,14.2)$.  Note that this point is in our shaded region,
but it is not in the region excluded by Miller's~\cite{Mil1} 
argument using the explicit formula.

In \eqref{eqn:beta12} we used the point $(\lambda_1,\lambda_2) = (8.4,14.2)$
to form one equation for the coefficients.  Now we will make another equation.
Choosing $g_3(s)= e^{is/4}$,  and continuing to use $s_0=\frac12+i$, we find
\begin{align} \label{eqn:beta3}
g_3(s_0)^{-1}  &\,\RHS(s_0,Q,\kappa,\lambda,\varepsilon,\{a_n\},g_3)= \\
&\mathstrut -1078.21 -2.38608\, b_1(2)-256.449 \, b_2(2)
-3.595\, b_1(3)+17.503 \, b_2(3) \cr
&+1.0549 \, b_1(4) -0.1226 \, b_2(4) 
+\cdots  -4.84 \cdot 10^{-6}\, b_1(8) +0.000012 \, b_2(8) + \cdots.\nonumber
\end{align}
Recall that $a(n)$ are the Dirichlet coefficients, and we write
$a(n)= b_1(n) + i b_2(n)$ with $b_1, b_2\in \R$.

Since \eqref{eqn:beta1},  \eqref{eqn:beta2}, and \eqref{eqn:beta3} 
each equal $Z(\frac12+ i)$, setting them all equal gives two linear
equations in the Dirichlet coefficients.  
We can solve for $b_1(2)$ and  $b_2(2)$ to obtain
\begin{align}
b_1(2) = \mathstrut & 0.1866 + 0.0137\,  b_1(3) + 0.2209\,  b_2(3) 
	+0.0122 \,  b_1(4) + 0.0380\,  b_2(4) +\cdots \cr
	&- 2.7 \cdot 10^{-7} \, b_1(8) - 2.5  \cdot 10^{-6}  \, b_2(8) + \cdots - 2.3 \cdot 10^{-13} \, b_1(14)
- 1.6  \cdot 10^{-12} \, b_2(14)  \cdots.\cr
b_2(2) = \mathstrut &
-4.1854 - 0.0170\, b_1(3) + 0.0629\, b_2(3) + 0.0034 \, b_1(4) -0.0015\, b_2(4) +\cdots \cr
&- 2.4  \cdot 10^{-8}\, b_1(8) + 1.4  \cdot 10^{-7}  \, b_2(8) + \cdots  
-1.3  \cdot 10^{-14}  \, b_1(14) + 6.1  \cdot 10^{-14} \, b_2(14)  \cdots. \cr
\end{align}
Now we impose the Ramanujan bound \eqref{eqn:preciseramanujan}
on all the coefficients on the right side of the
above equation, giving
\begin{align}
b_1(2) = \mathstrut & \phantom{\mathstrut -\mathstrut}0.1866 \pm 1.33 \cr
b_2(2) = \mathstrut & -4.1854 \pm 0.301 .
\end{align}
Note that the $\pm$ terms  are rigorous bounds, not error estimates.
Therefore \hbox{$b_2(2) < -3.8842$}, provided all the other Dirichlet
coefficients satisfy the Ramanujan bound.  Since the Ramanujan bound
$|a(2)|\le 3$ implies $|b_2(2)|\le 3$, 
we have a contradiction.  The conclusion is that if $(\lambda_1,\lambda_2)$
are the parameters in the $\Gamma$-factors of an $L-$function, then not all
the Dirichlet coefficients satisfy the Ramanujan bound.  Thus, 
assuming the Ramanujan bound, there is no
$L-$function with $(\lambda_1,\lambda_2) = (8.4,14.2)$.

We eliminated $(\lambda_1,\lambda_2) = (8.4,14.2)$ by solving only
two linear equations.  
Identical calculations over a closely spaced grid of points give the
shaded region in Figure~\ref{fig:SL3Z}.  
This proves Theorem~\ref{thm:excluded}.

It is possible to eliminate a larger
region by solving a larger system, but we have not attempted to explore
the limit of this method.

\bibliographystyle{plain}

\end{document}